**Extremal estimates for strongly additive and strongly multiplicative arithmetic functions**

Victor Volfson

ABSTRACT The paper considers estimates for some sums and products of functions of prime numbers. Several assertions on this topic have been proven. We also study extremal estimates for strongly additive and strongly multiplicative arithmetic functions. Several assertions on this topic are proved and examples are considered.

## 1. INTRODUCTION

An arithmetic function is a function defined on the set of natural numbers and taking values on the set of complex numbers. The name arithmetic function is due to the fact that this function expresses some arithmetic property of the natural series.

Estimating asymptotics of arithmetic functions and moments of arithmetic functions has been and remains at present an urgent problem [1,…,3].

Arithmetic functions usually have the specific property of changing internally irregularly and chaotically, as a result of which classical methods of analysis, as a rule, are powerless to adequately describe their behavior.

Probabilistic number theory responds to the naturally occurring task under such circumstances of making a statistical study of this behavior. The concept of the normal distribution of an arithmetic function is intended to reflect its "almost" certain behavior. In practice, this involves the process of neglecting a set of zero-density integers (obviously depending on the function in question) to eliminate erroneous values.

This "almost everywhere" approach was used in [4] to determine the asymptotics of additive and multiplicative arithmetic functions.

The main disadvantage of this approach is that asymptotics of arithmetic functions are determined with a probability equal to 1, and not exactly.

An approach from the point of view of extremal estimates of arithmetic functions will be considered in the paper, which is free from this drawback. This approach gives accurate estimates.

Let's look at some properties of arithmetic functions.

Let a natural number $m$ have a canonical decomposition $m = p_1^{a_1}...p_t^{a_t}$, where $p_i$ is a prime number and $\alpha_i$ is a natural number.

Then, the property for the additive arithmetic function is fulfilled:
$$f(m) = f(p_1^{a_1}...p_t^{a_t}) = f(p_1^{a_1}) + ... + f(p_t^{a_t}) = \sum_{p^\alpha | m} f(p^\alpha).$$

The following property holds for the corresponding strongly additive function:
$$f^*(m) = f^*(p_1^{a_1}...p_t^{a_t}) = \sum_{p|m} f(p).$$

The property for the multiplicative arithmetic function is executed:
$$g(m) = g(p_1^{a_1}...p_t^{a_t}) = \prod_{p^\alpha \| m} g(p^\alpha).$$

The following property holds for the corresponding strongly multiplicative function:
$$g^*(m) = g^*(p_1^{a_1}...p_t^{a_t}) = \prod_{p|m} g(p).$$

A multiplicative arithmetic function $g(m) > 0$ is converted to an additive one using the logarithm operation:
$$\ln(g(m)) = \ln(g(p_1^{a_1}...p_t^{a_t})) = \ln(\prod_{p^\alpha \| m} g(p)) = \sum_{p^\alpha \| m} \ln(g(p)).$$

An additive arithmetic function is converted into a multiplicative one using the potentiation operation: $e^{f(m)} = e^{f(p_1^{a_1}...p_t^{a_t})} = e^{f(p_1^{a_1})+...+f(p_t^{a_t})} = \prod_{p^\alpha \| m} e^{f(p^\alpha)}.$

Logarithm and potentiation can be carried out in any base $a > 0, a \neq 1$. We will consider more often $a = e$ in this work.

Based on these properties, extremal estimates for strongly additive and strongly multiplicative arithmetic functions are determined by estimates of certain sums and products of prime functions, which we will consider in the next chapter of this paper.

2. ESTIMATES FOR CERTAIN SUMS AND PRODUCTS OF FUNCTIONS OF PRIMES

Asymptotic estimates for sums of functions of prime numbers of the form $\sum_{p \leq n} f(p)$ were considered in [5],[6]. However, an extremal estimate of strongly additive arithmetic functions requires consideration of sums of the form $\sum_{p|n} f(p)$, where $f$ is the corresponding additive arithmetic function.

Estimating these sums in the form: $\sum_{p|n} f(p) \leq \sum_{p \leq n} f(p)$ is a rough estimate. Therefore, it makes sense to search for more accurate extremal estimates for the sums of functions of prime numbers of the form $\sum_{p|n} f(p)$. Moreover, the requirement that $f$ is an additive arithmetic function is optional.

Assertion 1

$$\sum_{p|n} f(p) \leq \omega(n) \sup_{2 \leq x \leq n}\{f(x)\}. \tag{2.1}$$

Proof

Let $n = \sum_{i=1}^{k} p_i^{\alpha_i}$ then:

$$\sum_{p|n} f(p) = f(p_1) + f(p_2) + ... + f(p_k) = \sum_{i=1}^{k} f(p_i). \tag{2.2}$$

Based on (2.2) and taking into account what $f(p_i) \leq \sup_{2 \leq x \leq n}\{f(x)\}$ we get:

$$\sum_{p|n} f(p) \leq k \sup_{2 \leq x \leq n}\{f(x)\}. \tag{2.3}$$

Since $k = \omega(n)$, then based on (2.3) we obtain (2.1).

Let's look at an example. Based on Assertion 1, we estimate:

$$\sum_{p|n} \ln \varphi(p) \leq \omega(n) \sup_{x \in [2,n]}\{\ln \varphi(x)\} \leq \omega(n) \ln(n) = O(\frac{\ln^2(n)}{\ln \ln(n)}). \tag{2.4}$$

Now let's do the evaluation in a different way:

$$\sum_{p|n} \ln \varphi(p) = \sum_{p|n} \ln(p-1) \leq \sum_{p \leq n} \ln(p-1) = \sum_{p \leq n} \ln p(1-1/p) = \sum_{p \leq n} \ln p + \sum_{p \leq n} \ln(1-1/p) = n + O(n/\ln(n)). \tag{2.5}$$

A comparison of estimates (2.4) and (2.5) shows that estimate (2.4) is more accurate.

Assertion 2

If an arithmetic function $f$ is monotonically decreasing, then

$$\sum_{p|n} f(p) \leqslant \sum_{p \leqslant p_{\omega(n)}} f(p), \tag{2.6}$$

where $p_k$ is the $k$-th prime number, and $p_{\omega(n)} \leqslant (1+o(1))\ln n$.

Proof

The estimate $\sum_{p|n} f(p) \leqslant \sum_{p \leqslant p_{\omega(n)}} f(p)$, follows from the fact that both sums have the same number of terms, but the terms in the right sum are not less (due to the monotonic decrease $f$).

The estimate $p_{\omega(n)} \leqslant (1+o(1))\ln n$ follows from the inequality:

$$\ln n \geqslant \ln \prod_{p \leqslant p_{\omega(n)}} p = \sum_{p \leqslant p_{\omega(n)}} \ln p = \theta(p_{\omega(n)}) = p_{\omega(n)}(1+o(1))$$

(where $\theta(x)$ is the Chebyshev theta function) and the asymptotic law of distribution of prime numbers in the form $\theta(x) \sim x$ when the value $x \to +\infty$.

Example for assertion 2. Let's evaluate the following arithmetic function:

$$\sum_{p|n} \ln \frac{n}{\varphi(n)} = \sum_{p|n} \ln \frac{p}{p-1} = \sum_{p|n} \ln(1+\frac{1}{p-1}) = \sum_{p|n}(\frac{1}{p-1} + O(\frac{1}{(p-1)^2})) = \sum_{p|n} \frac{1}{p-1} + O(1). \quad (2.7)$$

Based on (2.7) and Assertion 2, we obtain:

$$\sum_{p|n} \ln \frac{n}{\varphi(n)} = \sum_{p|n} \frac{1}{p-1} + O(1) \leq \sum_{p \leq \ln(n)(1+o(1))} \frac{1}{p-1} + O(1) = \ln\ln\ln(n)(1+o(1)) + O(1). \quad (2.8)$$

Assertion 3

Let the arithmetic function $g(m) > 0, m = 1,...,n$ then

$$\prod_{p|n} g(p) \leq e^{\omega(n)\sup_{x \in [2,n]}\{\ln g(x)\}}. \quad (2.9)$$

Proof

We take the logarithm of the product of arithmetic functions:

$$\ln(\prod_{p|n} g(p)) = \sum_{p|n} \ln g(p). \quad (2.10)$$

Let's apply assertion 1 to (2.10):

$$\sum_{p|n} \ln g(p) \leq \omega(n)\sup_{x \in [2,n]}\{\ln g(x)\}. \quad (2.11)$$

Using (2.11) we get:

$$e^{\sum_{p|n} \ln g(p)} = \prod_{p|n} g(p) \leq e^{\omega(n)\sup_{x\in[2,n]}\{\ln g(x)\}},$$

which corresponds to (2.9).

Now we will consider an example for assertion 3. Let's make an estimate of the product of arithmetic functions $\varphi(m) > 0, m = 1,...,n$:

$$\prod_{p|n} \varphi(p) \leq e^{\omega(n)\sup_{x\in[2,n]}\{\ln \varphi(x)\}} \leq e^{\omega(n)\ln(n)} = n^{\omega(n)} \leq n^{c\ln(n)/\ln\ln(n)}, \qquad (2.12)$$

where c is a constant.

Assertion 4

Let the arithmetic function $g(m) > 0, m = 1,...,n$ and $A(n) = e^{\sum_{p \leq p_w(n)} \ln(g(p))}$, where $p_\omega(n) \leq (1+o(1))\ln n$, then, if $g(m)$ decreases monotonically, then the following upper bound is fulfilled:

$$\prod_{p|n} g(p) \leq A(n), \qquad (2.13)$$

Proof

We take the logarithm of the product of arithmetic functions:

$$\ln \prod_{p|n} g(p) = \sum_{p|n} \ln g(p). \qquad (2.14)$$

If $g(m)$ decreases monotonically, then $\ln g(m)$ also decreases monotonically, therefore, based on Assertion 2 and (2.14), we obtain:

$$\ln \prod_{p|n} g(p) = \sum_{p|n} \ln g(p) \leq \sum_{p \leq p_\omega} \ln g(p). \qquad (2.15)$$

Based on (2.15) we have:

$$e^{\ln \prod_{p|n} g(p)} = \prod_{p|n} g(p) \leq e^{\sum_{p \leq p_\omega} \ln g(p)},$$

which corresponds to (2.13).

Пример на утверждение 4.

Let's look at an example for assertion 4.

The arithmetic function $g(p) = \dfrac{p}{p-1} = 1 + \dfrac{1}{p-1}$ is positive and monotonically decreasing, so Assertion 4 can be applied to it:

$$\prod_{p|n} \frac{p}{p-1} \leq e^{\sum_{p \leq p_\omega} \ln \frac{p}{p-1}}. \qquad (2.16)$$

Based on (2.7), (2.8) and (2.16) we get:

$$\frac{n}{\varphi(n)} = \prod_{p|n} \frac{p}{p-1} \leq e^{\sum_{p \leq p_\omega} \ln \frac{p}{p-1}} = e^{\ln\ln\ln n(1+o(1)) + O(1)} \leq e^c \ln\ln n (1+o(1)), \qquad (2.17)$$

where c is a constant.

Having in mind (2.17), we obtain a lower estimate for the Euler function:

$$\varphi(n) \geq \frac{e^{-c} n (1+o(1))}{\ln\ln n}. \qquad (2.18)$$

Now let's talk about lower bounds for some sums and products of functions of primes.

Assertion 5

Let the function $f(p) > 0$ monotonically increase at $p \geq a$, where $p, a$ are prime numbers, then:

$$\sum_{p|n} f(p) \geq f(a). \qquad (2.19)$$

Equality is achieved when the value is $n = a^m$, where $m$ is a natural number.

Proof

The minimum value $\sum_{p|n} f(p)$ is reached at a value $n$ equal to the power of a prime number, since the sum increases for several different prime divisors of the number $n$, having in mind, that $f(p) > 0$.

Since the minimum prime number from which the function $f(p)$ monotonically increases is $a$, then the minimum value $\sum_{p|n} f(p)$ is reached at the value $n = a^m$, where $m$ is a natural number. Therefore, inequality (2.19) holds.

Consequence of Assertion 5

Let the functions $f_i(p) > 0$ and monotonically increase as $p \geq a_i$, where $p, a_i$ are prime numbers, and $i = \{1, 2, ..., k\}$. Then the following inequality holds:

$$\sum_{i \leq k} \sum_{p|n} f_i(p) \geq \sum_{i \leq k} f_i(a),$$

where $a = \sup\{a_1, a_2, ..., a_k\}$.

Proof

Since $f_i(p) > 0$ and monotonically increase as $p \geq a_i$, then $\sum_{i \leq k} f_i(p) > 0$ and also monotonically increases as $p \geq \sup\{a_1, a_2, ..., a_k\}$ and, therefore, satisfies the conditions of Assertion 5.

Example for assertion 5.

The function $f(p) = \dfrac{p^2 + p + 1}{p - 1} > 0$ is monotonically increasing at $p \geq 3$, so the following lower estimate is made for the following sum:

$$\sum_{p|n} f(p) \geq f(3) = 6{,}5. \qquad (2.20)$$

Equality in (2.20) is achieved at the value $n = 3^m$.

Assertion 6

Let the function $g(p) > 1$ monotonically increase at $p \geq a$, where $a$ is a prime number, then:

$$\prod_{p|n} g(p) \geq g(a). \qquad (2.21)$$

Equality is achieved when the value is $n = a^m$, where $m$ is a natural number.

Proof

Let's take the logarithm $\prod_{p|n} g(p)$ and get:

$$\ln \prod_{p|n} g(p) = \sum_{p|n} \ln g(p). \qquad (2.22)$$

Since the function $g(p) > 1$ and monotonically increases as $p \geq a$, then $\ln g(p) > 0$ and it also monotonically increases as $p \geq a$, therefore, the function $\ln g(p)$ satisfies all the conditions of Assertion 5. Therefore, based on (2.22), the lower bound is fulfilled:

$$\sum_{p|n} \ln g(p) \geq \ln g(a). \qquad (2.23)$$

Taking into account (2.23), we can write:

$$e^{\sum_{p|n} g(p)} = \prod_{p|n} g(p) \geq e^{\ln g(a)} = g(a),$$

which corresponds to (2.21).

Consequence from Assertion 6

Let functions $g_i(p) > 1$ and monotonically increase at $p \geq a_i$, where $p, a_i$ are prime numbers, and $i = \{1, 2, ..., k\}$, then the following inequality is true:

$$\prod_{i \leq k} \prod_{p|n} g_i(p) \geq \prod_{i \leq k} g_i(a),$$

where $a = \sup\{a_1, a_2..., a_k\}$.

Proof

Since functions $g_i(p) \geq 1$ monotonically increase as $p \geq a_i$, then $\prod_{i \leq k} g_i(p) > 1$ monotonically increases as $p \geq a$, where $a = \sup\{a_1, a_2..., a_k\}$ and satisfies the conditions of Assertion 6.

Example for Assertion 6.

The function $g(p) = \dfrac{p^2 + p + 1}{p - 2} > 1$ monotonically increases as $p \geq 5$, so the following lower bound is fulfilled:

$$\prod_{p|n} g(p) \geq g(5) = \frac{31}{3}. \tag{2.24}$$

Equality in (2.24) is achieved at the value $n = 5^m$.

## 3. EXTREMAL ESTIMATES FOR STRONGLY ADDITIVE AND MULTIPLICATIVE ARITHMETIC FUNCTIONS

Assertion 7

Let $f(m), m = 1,...,n$ is an additive arithmetic function and if $f(p^\alpha) \geq f(p)$, then it is performed for a strictly additive arithmetic function,:

$$\sum_{p|n} f(p) \leq f(n), \tag{3.1}$$

if $f(p^\alpha) \leq f(p)$, then:

$$\sum_{p|n} f(p) \geq f(n). \tag{3.2}$$

Proof

Let $f(n) = f(p_1^{\alpha_1} p_2^{\alpha_2} ... p_k^{\alpha_k})$, then, by virtue of additivity $f(m), m = 1,...,n$, we obtain:

$$f(n) = f(p_1^{\alpha_1} p_2^{\alpha_2} ... p_k^{\alpha_k}) = f(p_1^{\alpha_1}) + f(p_2^{\alpha_2}) + ... + f(p_k^{\alpha_k}) = \sum_{p|n} f(p^\alpha). \tag{3.3}$$

If $f(p^\alpha) \geq f(p)$, then by virtue of (3.3):

$$f(n) = \sum_{p|n} f(p^\alpha) \geq \sum_{p|n} f(p) \text{ or } \sum_{p|n} f(p) \leq f(n),$$

which corresponds to (3.1).

If $f(p^\alpha) \leq f(p)$, then by virtue of (3.3):

$$f(n) = \sum_{p|n} f(p^\alpha) \leq \sum_{p|n} f(p) \text{ or } \sum_{p|n} f(p) \geq f(n),$$

which corresponds to (3.2).

Example for Assertion 7

The arithmetic function $\ln \varphi(n)$ is an additive arithmetic function and it is executed for it when the value is $\alpha \geq 2$:

$$\ln \varphi(p^\alpha) = \ln p^{\alpha-1}(p-1) \geq \ln(p-1) = \ln \varphi(p).$$

Therefore, based on Assertion 7 (3.1), we obtain the following estimate:

$$\sum_{p|n} \ln \varphi(p) \leq \ln \varphi(n) \leq \ln n. \tag{3.4}$$

Estimate (3.4) is more accurate than estimate (2.4) and even more so than estimate (2.5).

The function $\ln \varphi(p) = \ln(p-1) \geq 0$ and monotonically increases as $p \geq 2$, therefore, based on Assertion 5, a lower bound is obtained $\sum_{p|n} \ln \varphi(p) \geq 0$ for a strongly additive arithmetic function $\sum_{p|n} \ln \varphi(p)$. Equality is achieved with the value $n = 2^m$.

One more example. Let's consider an arithmetic function $\ln \frac{\tau(n)}{n}$ with the value $n \geq 2$. The function $\ln \frac{\tau(n)}{n}$ is an additive arithmetic function and it is executed when the value is $n \geq 2$:

$$\ln \frac{\tau(p^\alpha)}{p^\alpha} = \ln \frac{\alpha+1}{p^\alpha} \leq \ln \frac{2}{p} = \ln \frac{\tau(p)}{p}.$$

Therefore, based on Assertion 7 (3.2) and the fact that $\tau(n) \geq 2$, the following holds:

$$\sum_{p|n} \ln \frac{\tau(p)}{p} \geq \ln \frac{\tau(n)}{n} \geq \ln \frac{2}{n}. \tag{3.5}$$

Equality is achieved with the value $n = p$.

On the other hand, the arithmetic function $\ln \frac{\tau(p)}{p} = \ln \frac{2}{p}$ is monotonically decreasing, so Assertion 2 can be applied:

$$\sum_{p|n} \ln \frac{\tau(p)}{p} \leq \sum_{p \leq \ln n(1+o(1))} \ln \frac{2}{p} = \ln 2 \sum_{p \leq p_{w(n)}} 1 - \sum_{p \leq p_{w(n)}} \ln p = \ln 2 \frac{p_{w(n)}}{\ln p_{w(n)}} - p_{w(n)} = -p_{w(n)}(1 - \frac{\ln 2}{\ln p_{w(n)}}), (3.6)$$

where $p_{w(n)} = \ln n(1+o(1))$.

Based on (3.5) and (3.6), we obtain an asymptotic estimate for the value $n \to \infty$:

$$\ln \frac{2}{n} \leq \sum_{p|n} \ln \frac{\tau(p)}{p} \leq \ln \frac{1}{n}(1 - \frac{\ln 2}{\ln \ln n}).$$

A lower bound for the sum of strongly additive functions can be carried out on the basis of the corollary of Assertion 5.

For example, the function $\ln \varphi(p) = \ln(p-1)$ is greater than or equal to 0 and monotonically increases for $p \geq 2$.

The function $\ln \sigma(p) = \ln(p+1)$ is greater than 0 and also increases monotonically for $p \geq 2$.

Therefore, based on the corollary of Assertion 5, we obtain a lower estimate for the sum of strongly additive arithmetic functions:

$$\sum_{p|n} \ln \varphi(p) + \sum_{p|n} \ln \sigma(p) = \prod_{p|n} \ln(p-1) + \prod_{p|n} \ln(p+1) \geq \ln 3.$$

Assertion 8

Let $g(m) > 0, m = 1,...,n$ is a multiplicative arithmetic function and if $g(p^\alpha) \leq g(p)$, then it is performed for a strictly multiplicative arithmetic function $\prod_{p|n} g(p)$:

$$\prod_{p|n} g(p) \leq g(n) \quad (3.7)$$

and vice versa, if $g(p^\alpha) \geq g(p)$, then

$$\prod_{p|n} g(p) \geq g(n). \tag{3.8}$$

Proof

Let $g(m) > 0, m = 1,...,n$ is a multiplicative arithmetic function, then $f(m) = \ln g(m), m = 1,...,n$ is the corresponding additive arithmetic function and vice versa $g(m) = e^{f(m)}, m = 1,...,n$.

Therefore $g(p^\alpha) = e^{f(p^\alpha)}$, $g(p) = e^{f(p)}$ and the condition $g(p^\alpha) \geq g(p)$ can be written as:

$$g(p^\alpha) = e^{f(p^\alpha)} \geq g(p) = e^{f(p)}. \tag{3.9}$$

It follows from (3.9) that $f(p^\alpha) \geq f(p)$, therefore, based on Assertion 5, the following holds:

$$\sum_{p|n} f(p) \leq f(n) \text{ or } e^{\sum_{p|n} f(p)} = \prod_{p|n} g(p) \leq e^{f(n)} = g(n),$$

which corresponds to (3.7).

The condition $g(p^\alpha) \leq g(p)$ can be written as:

$$g(p^\alpha) = e^{f(p^\alpha)} \leq g(p) = e^{f(p)}. \tag{3.10}$$

It follows from (3.10) that $f(p^\alpha) \leq f(p)$, therefore, based on Assertion 5, we have:

$$\sum_{p|n} f(p) \geq f(n) \text{ or } e^{\sum_{p|n} f(p)} = \prod_{p|n} g(p) \geq e^{f(n)} = g(n),$$

which corresponds to (3.8).

Example for assertion 8.

Let's consider a multiplicative function $\sigma_k(n) = \sum_{d|n} d^k$ - the sum of the divisors of the $k$ th power of a natural number. Note that when $\alpha \geq 2$:

$$\sigma_k(p^\alpha) = \frac{p^{(\alpha+1)k} - 1}{p^k - 1} \geq \frac{p^{2k} - 1}{p^k - 1} = p^k + 1 = \sigma_k(p). \tag{3.11}$$

Based on Assertion 8 and (3.11) it follows:

$$\prod_{p|n} \sigma_k(p) = \prod_{p|n}(p^k+1) \leq \sigma_k(n). \tag{3.12}$$

Taking into account the maximum order $\sigma_k(n)$ [7], based on (3.12), we obtain:

$$\prod_{p|n} \sigma_k(p) = \prod_{p|n}(p^k+1) \leq \sigma_k(n) \leq n^k \varsigma(k), k > 1, \tag{3.13}$$

where $\varsigma(k)$ is the value of the Riemann zeta function at natural value $k$.

Based on (3.12) and taking into account the maximum order [7] at value $k=1$ we obtain:

$$\prod_{p|n} \sigma(p) = \prod_{p|n}(p+1) \leq \sigma(n) \leq e^\gamma n \ln\ln(n). \tag{3.14}$$

Based on (3.12) and taking into account the maximum order [7] at value $k=0$ we obtain:

$$\prod_{p|n} \tau(p) = \prod_{p|n} 2 \leq \tau(n) \leq n^{\frac{\ln 2}{\ln\ln(n)}}. \tag{3.15}$$

We use assertion 4 to get the estimate $\prod_{p|n} \tau(p)$:

$$\prod_{p|n} \tau(p) = \prod_{p|n} 2 \leq e^{\sum_{p \leq p_{w(n)}} \ln 2} = e^{\ln 2 \frac{p_{w(n)}}{\ln p_{w(n)}}}, \tag{3.16}$$

where $p_{w(n)} = \ln n(1+o(1))$.

The resulting asymptotic estimate (3.16) corresponds to (3.15) when the value $n \to \infty$.

A lower estimate is known for the multiplicative arithmetic function $\sigma_k(n) \geq n^k, k > 0$.

However, the lower bound for the corresponding strongly multiplicative arithmetic function $\prod_{p|n} \sigma_k(p) = \prod_{p|n}(p^k+1)$ is significantly less than $n^k, k > 0$, at points with prime powers greater than 1 (not square-free).

Let us use the fact that the function $\sigma_k(p) = p^k + 1$ is greater than or equal to 0 and monotonically increases as $p \geq 2$.

Based on Assertion 6, the lower bound for a strongly multiplicative arithmetic function $\prod_{p|n} \sigma_k(p) = \prod_{p|n}(p^k + 1)$ is:

$$\prod_{p|n} \sigma_k(p) = \prod_{p|n}(p^k + 1) \geq 2^k + 1. \tag{3.17}$$

The minimum value of a strongly multiplicative arithmetic function $\prod_{p|n} \sigma_k(p) = \prod_{p|n}(p^k + 1)$ is reached when $n = 2^m$ and is equal to $2^k + 1$.

A lower estimate for the product of strongly multiplicative arithmetic functions can be carried out on the basis of the corollary of Assertion 6.

For example, the function $\varphi(p) = p - 1$ is greater than 1 and monotonically increases as $p \geq 2$.

The function $\sigma_2(p) = p^2 + 1$ is greater than 1 and monotonically increases as $p \geq 2$.

Based on the corollary of Assertion 6, we obtain a lower bound for the product of strongly multiplicative functions:

$$\prod_{p|n} \varphi(p) \prod_{p|n} \sigma_2(p) \geq \varphi(2)\sigma_2(2) = 5.$$

## 4. CONCLUSION AND SUGGESTIONS FOR FURTHER WORK

The next article will continue to study the asymptotic behavior of some arithmetic functions.

## 5. ACKNOWLEDGEMENTS

Thanks to everyone who has contributed to the discussion of this paper. I am grateful to everyone who expressed their suggestions and comments in the course of this work.

References

1. J.Chattopadhyay, P. Darbar Mean values and moments of arithmetic functions over number fields, Research in Number Theory 5(3), 2019.

2. M. Garaev, M. Kühleitner, F. Luca, G. Nowak Asymptotic formulas for certain arithmetic functions, Mathematica Slovaca 58(3):301-308, 2008.

3. R. Brad Arithmetic functions in short intervals and the symmetric group, arXiv preprint https://arxiv.org/abs/1609.02967(2016).

4. Victor Volfson Asymptotic behavior of additive and multiplicative arithmetic functions, https://arxiv.org/abs/2209.06907(2022).

5. Victor Volfson Asymptotics for sums of functions of primes, arXiv preprint https://arxiv.org/abs/2008.11100 (2020).

6. Victor Volfson Asymptotics of probability characteristics of additive arithmetic functions, arXiv preprint https://arxiv.org/abs/2108.13199 (2021).

7. Gérald Tenenbaum. Introduction to analytic and probabilistic number theory. Third Edition, Translated by Patrick D. F. Ion, 2015 - 630 p.